\documentclass[a4paper, reqno, 10pt]{amsart}
\usepackage[pdftex,pagebackref,letterpaper=true,colorlinks=true,pdfpagemode=none,urlcolor=blue,linkcolor=blue,citecolor=blue,pdfstartview=FitH]{hyperref}

\usepackage{amsmath,amssymb}
\usepackage{graphicx}
\usepackage{cite}
\usepackage[usenames,dvipsnames]{xcolor}
\newcommand{\blue}{\textcolor{blue}}

\newtheorem{theorem}{Theorem}

\newtheorem{corollary}[theorem]{Corollary}

\newtheorem{conjecture}[theorem]{Conjecture}

\begin{document}

  \title[A Friendly Intro to Sieves]{A Friendly Intro to Sieves with a Look Towards Recent Progress on the Twin Primes Conjecture}

  \author{David Lowry-Duda}
  \address{ Brown University Mathematics}
  \email{djlowry@math.brown.edu}
  \thanks{This material is based upon work supported by the National Science Foundation Graduate Research Fellowship Program under Grant No. DGE-0228243}


  \date{\today}

  \maketitle

  This is an extension and background to a talk I gave on 9 October 2013 to the Brown Graduate Student Seminar, called `A friendly intro to sieves with a look towards recent progress on the twin primes conjecture.' During the talk, I mention several sieves, some with a lot of detail and some with very little detail. I also discuss several results and built upon many sources. I'll provide missing details and/or sources for additional reading here. 

  Furthermore, I like this talk, so I think it's worth preserving.

  \section{Introduction}

  We talk about sieves and primes. Long, long ago, Euclid famously proved the infinitude of primes ($\approx 300$ B.C.). Although he didn't show it, the stronger statement that the sum of the reciprocals of the primes diverges is true:
  \[
    \sum_{p} \frac{1}{p} \to \infty,
  \]
  where the sum is over primes.

  \begin{proof}
    Suppose that the sum converged. Then there is some $k$ such that 
    \[
      \sum_{i = k+1}^\infty \frac{1}{p_i} < \frac{1}{2}.
    \]
    Suppose that $Q := \prod_{i = 1}^k p_i$ is the product of the primes up to $p_k$. Then the integers $1 + Qn$ are relatively prime to the primes in $Q$, and so are only made up of the primes $p_{k+1}, \ldots$. This means that
    \[
      \sum_{n = 1}^\infty \frac{1}{1+Qn} \leq \sum_{t \geq 0} \left( \sum_{i > k} \frac{1}{p_i} \right) ^t < 2,
    \]
    where the first inequality is true since all the terms on the left appear in the middle (think prime factorizations and the distributive law), and the second inequality is true because it's bounded by the geometric series with ratio $1/2$. But by either the ratio test or by limit comparison, the sum on the left diverges (aha! Something for my math 100 students), and so we arrive at a contradiction.

    Thus the sum of the reciprocals of the primes diverges.
  \end{proof}

  I learned of this proof from Apostol's Introduction to Analytic Number Theory, and although I've seen many proofs since, this is still my favorite.

  This turned out to be a pretty good proof technique. In the 1800's, Dirichlet and others proved the Prime Number Theorem, and the more sophisticated prime number theorem for primes in arithmetic progressions (sometimes called Dirichlet's Theorem). The first says that the number of primes less than $x$, which I'll denote by $\pi(x)$, is asymptotically $\frac{x}{\log x}$. The second says that as long as $\gcd(a,b) = 1$, then the arithmetic progression $a, a+b, a+2b, a+3b, \ldots$ contains infinitely many primes.

  Sort of similar to before, a stronger statement is true:

  \[
    \sum_{a + bk = p} \frac{1}{a + bk} \to \infty,
  \]
  where the summation is just over those elements where $a + bk$ is prime. Further, Dirichlet showed that each progressions mod $a$ have the same asymptotic, and so primes are very equidistributed. 

  \begin{quote}
    As an aside, the primes are \emph{very equidistributed} in the sense that the decimal $0.23571113\ldots$, gotten from concatenating all the primes, is a \textbf{normal number}, meaning that every finite pattern of digits appears in the decimal, and every pattern of the same length occurs with essentially the same density. So no digit, digits, or pattern of digits appear any more often than any other digit, digits, or pattern of digits. This number is called the Copeland-Erd\"{o}s constant.
  \end{quote}

  A few decades later, in 1859, Riemann wrote and published his famous Memoir. In this, he introduced what we now call the Riemann zeta function $\zeta(s) = \sum_{n \geq 1} \frac{1}{n^s}$, gave its analytic continuation and functional equation, and created the field of analytic number theory. This matters to me, because I am an analytic number theorist, and it's good to know your roots.

  Perhaps these methods and developments are what inspired Viggo Brun to try to analyze twin primes around the start of the 20th century. Twin primes are pairs of primes that are $2$ apart. For instance, $3$ and $5$ are twin primes, as are $11$ and $13$, and so on. We know there are infinitely many primes - are there infinitely many twin primes? What about ``cousin'' primes (pairs of primes of the form $p, p+4$) or ``sexy primes'' ($p, p+6$)? Or what about bigger sets of primes. Are there infinitely many trios of primes $p, p+2, p+4$, or $p, p+2, p+6$?

  \begin{conjecture}[Twin Primes Conjecture]
    There are infinitely many pairs of twin primes $p, p+2$.
  \end{conjecture}

  Brun wanted to analyze the sum
  \[
    \sum_{p \text{  twin}} \frac{1}{p},
  \]
  perhaps hoping that the sum would be infinite and thus showing that there are infinitely many twin primes. Brun successfully analyzed this sum, but he did not manage to prove or disprove the Twin Primes Conjecture. Instead, he showed that
  \[
    \sum_{p \text{  twin}} \frac{1}{p} \approx 1.9 < \infty.
  \]
  This sum is \emph{finite}, so if there are infinitely many twin primes, then there aren't too too many of them. On the other hand, just because the sum is finite doesn't mean there are only finitely many twin primes. For example,
  \[
    \sum_{n \geq 1} \frac{1}{n^2} = \frac{\pi^2}{6},
  \]
  a finite number, and there are clearly infinitely many squares. (Showing this equality is true is called Basel's problem and it is a classic problem in a complex analysis class).

  Brun's result is impressive, but it's not sufficient to say anything about the infinitude of twin primes. For a long time, it was widely thought that no one was getting any closer to proving something about the infinitude of twin primes than Brun was over a century ago. But then in 2013, Yitang Zhang broke the stalemate by showing that there are infinitely many primes of the form $p, p+2k$ for some fixed $k$ (although he didn't prove what that $k$ was). Shortly afterwards, James Maynard showed that there are infinitely primes of the form $p, p + a, p + b$ for some fixed $a, b$ (in fact, he proved a much stronger result stating that there are infinitely many prime families of many types). 

  The most up-to-date results and progress is being carried out by Terry Tao and the \href{http://michaelnielsen.org/polymath1/index.php?title=Bounded_gaps_between_primes}{Polymath8} massively collaborative mathematics project.

  \section{The Sieve of Eratosthenes}

  Brun, Zhang, Maynard, and Polymath8 all worked with \emph{sieves}. In some ways, math sieves are just like sieves you might see in a kitchen: they filter some stuff out, and hopefully let what you want through. They don't work perfectly. Some extra stuff usually gets through, or you don't get everything you want, or (most likely) you get a little bit of both.

  Many have heard of a sieve used in mathematics. Around 250 BCE, Greek polymath and librarian Eratosthenes of Cyrene developed a sieve to find and count primes. In his honor and memory, we call it the Sieve of Eratosthenes. To understand his sieve, let us first try to find and count primes ourselves.

  The naive method would be to use ``trial division'' on each number to check if it's composite. If not, then it's a prime. Then we go to the next number. For example, we might wonder if $57$ is prime. Is it divisible by $2$? No. Is it divisible by $3$? Yes! Oh - so it's not prime. Then we might check $58$, then $59$, and so on. 

  This can be improved by checking only if a number $n$ is divisible by primes $p \leq \sqrt n$, since any composite number has at least one prime factor less than its square root. This would save time checking $59$, for example, since you would check $2, 3, 5,$ and $7$, and since it's not divisible by any of those, we know it's prime. 

  We can save a little more time by noticing that we can skip every even number after $2$, since they must be divisible by $2$. If we think about it, we see we can skip every multiple of $3$ after $3$, too, for the same reason. They must be divisible by $3$. And every multiple of $4$ after $4$ - but since $4$ is a multiple of $2$, this step is redundant.

  These are the insights that led to the Sieve of Eratosthenes. To find the primes up to $25$, we first write down the potential numbers, 
  \[
    1, 2, 3, 4, 5, 6, 7, 8, 9, 10, 11, 12, 13, 14, 15, 16, 17, 18, 19, 20, 21, 22, 23, 24, 25.
  \]
  Since $1$ isn't prime, cross it out. ($1$ is going to be a pain today, as it doesn't quite fit all the patterns)
  \[
    \not 1, 2, 3, 4, 5, 6, 7, 8, 9, 10, 11, 12, 13, 14, 15, 16, 17, 18, 19, 20, 21, 22, 23, 24, 25.
  \]
  The next number on the list will be prime. So $2$ is prime. We'll underline it. Then cross out all multiples of $2$.
  \[
    \not 1, \underline{2}, 3, \not 4, 5, \not 6, 7, \not 8, 9, \not 10, 11, \not 12, 13, \not 14, 15, \not 16, 17, \not 18, 19, \not 20, 21, \not 22, 23, \not 24, 25.
  \]
  Now we repeat. The next number on the list will be prime, and it's multiples should be crossed out.
  \[
    \not 1, \underline{2}, \underline{3}, \not 4, 5, \not 6, 7, \not 8, \not 9, \not 10, 11, \not 12, 13, \not 14, \not 15, \not 16, 17, \not 18, 19, \not 20, \not 21, \not 22, 23, \not 24, 25.
  \]
  Since $\sqrt{25} = 5$, and the next element on the list is $5$, this is our last step. We underline $5$, cross out any multiples of $5$ that are left, and everything left must be a prime. So we get
  \[
    \not 1, \underline{2}, \underline{3}, \not 4, \underline{5}, \not 6, \underline{7}, \not 8, \not 9, \not 10, \underline{11}, \not 12, \underline{13}, \not 14, \not 15, \not 16, \underline{17}, \not 18, \underline{19}, \not 20, \not 21, \not 22, \underline{23}, \not 24, \not 25.
  \]
  The primes are $2, 3, 5, 7, 11, 13, 17, 19, 23$. In particular, there are $9$ primes here, and $6$ primes bigger than $5$ (the largest prime we used to find them).

  This is must faster than trial division. And this is the plan of the Sieve of Eratosthenes. A different way of viewing the Sieve of Eratosthenes is that if we know the primes up to $5$ (which we do), then then we find the primes up to $5^2 = 25$ very quickly. If we wanted, we could repeat the process, using the primes up to $25$ to get the primes up to $25^2$ very quickly, and so on. So we can find the primes between $\sqrt{n}$ and $n$ very quickly.

  This is a form of the Sieve of Eratosthenes that gives the primes explicitly, but what if we were just interested in counting the number of primes (just as we wonder if there are infinitely many twin primes, but don't necessarily care to find all of them)? Let's look at our argument again with the Principle of Inclusion and Exclusion in mind. 

  Let $\pi(x, z)$ denote number of integers less than or equal to $x$ that are coprime to primes less than or equal to $z$, or rather

  \[
    \pi(x, z) = \#\left\{ n \leq x:p \not | \; n \; \forall p \leq z \right\} .
  \]

  So we want to try to figure out $\pi(n, \sqrt n)$ using the Sieve of Eratosthenes. First, let's look at our example to find $\pi(25, 5)$. We start with $25$ numbers. We first get rid of multiples of $2$. How many multiples of $2$ are there? There are $\lfloor 25 / 2 \rfloor = 12$ multiples of $2$ here (we are counting $2$ itself!). The next prime to remove is $3$. How many multiples of $3$ are there? There are $\lfloor 25/3 \rfloor = 8$ multiples of $3$. But wait, we've double counted a few things. 

  For example, we've counted the number $6$ twice, since it is both a multiple of $2$ and $3$. To not overcount, let's put back in those multiples of both $2$ and $3$, or rather let's add in the multiples of $6$. How many multiples of $6$ are there? There are $\lfloor 25/6 \rfloor = 4$. Then we take out the multiples of $5$: $\lfloor 25 / 5 \rfloor = 5$. But we've again overcounted, now by multiples of $10 = 2 \cdot 5$ and multiples of $15 = 3 \cdot 5$. So we add back in $\lfloor 25 / 10 \rfloor$ and $\lfloor 25 / 15 \rfloor = 1$. In principle, we would need to add back in multiples of $2 \cdot 3 \cdot5 = 30$, but since there are none less than $25$, that's not necessary.

  By the idea of the Sieve of Eratosthenes, this should be enough to eliminate all the composites between $5$ and $25$, and since we've also eliminated the primes up to $5$, we should be left with $\pi(25, 5)$.

  All together, this gives us
  \[
    25 - 12 - 8 - 5 + 4 + 2 + 1 - 0 = 7,
  \]
  so we expect that $\pi(25, 5) = 7$. But it happens to be that we counted the number of primes greater than $5$ and less than $25$ above, and we got only $6$. Why do we get $7$ instead? Well, it's an annoying thing: we never took out $1$ from the list of numbers when we were calculating $\pi(25, 5)$. On the one hand, this is a bit silly. On the other hand, $1$ is technically coprime to all primes less than $5$, and so should be included in this phrasing of the definition. Regardless, it turns out that it makes formulas easier to write down and consider if we include $1$ in this count, so that $1$ is almost a ``prime'' today.

  Number theorists have a function that makes writing this inclusion-exclusion counting argument easier. It's called the M\"{o}bius function $\mu(n)$, which is given by
  \[
    \mu(n) = \begin{cases}
      1 & n = 1 \\
      (-1)^k = (-1)^{\nu(n)} & n = p_1 p_2 \ldots p_k \\
      0 & p^2 \mid n \text{ for any } p
    \end{cases}.
  \]
  In other words, $\mu(n)$ is $0$ if any prime divides $n$ twice, and otherwise if $(-1)$ raised to the number of different primes dividing $n$. $\nu(n)$ happens to be another number theory shorthand, and it stands for the number of different primes dividing $n$. 

  So $\nu(p) = -1$ for any prime $p$, $\nu(p_1 p_2) = 1$ for any two distinct primes $p_1, p_2$, and so on. One final piece of notation: let $P(z) = \prod_{p \leq z} p$ be the product of the primes up to $z$.

  Then our Sieve of Eratosthenes style argument above can be written succinctly as
  \[
    \pi(25, 5) = \sum_{d \mid P(5)} \mu(d) \left\lfloor \frac{25}{d} \right\rfloor.
  \] 
  \emph{(If this is your first time seeing some of this notation, or if you don't believe it, I encourage you to write this expression out cleanly. You'll see that it's the same as the expression we have above giving $7$)}. This argument easily generalizes, so that

  \begin{theorem}[Original \emph{Counting} Sieve of Eratosthenes]
    Suppose $\pi(x,z)$ is the number of integers up to $x$ that are not divisible by primes up to $z$. Let $P(z)$ is the product of the primes up to $z$. Then we have
    \begin{equation}
      \pi(x, z) = \sum_{d \mid P(z)} \mu(d) \left\lfloor \frac{x}{d} \right\rfloor.
    \end{equation}
  \end{theorem}

  Clearly this is not the same form in which Eratosthenes would have presented this result, but the heart of it is as it was millennia ago. We are using multiplicative properties of integers (divisibility and properties related to the M\"{o}bius function in this case) to sift out particular numbers (smaller primes and all composites) and to isolate a set (larger primes) that we are interested in.

  \section{More General Sieves}

  Let us do as mathematicians tend to do and generalize. Let me be the first to say that what we will do is a bit technical and not obvious. But the key intuition is the same as in the Sieve of Eratosthenes. One great thing about the Sieve of Eratosthenes is that there is no error - it is an exact counting tool. But for more interesting or elusive subsets of the integers, we won't be able to be so precise. By giving up exactness, we will be able to estimate the sizes of a bigger family of subsets of the integers and not just large primes.

  In what follows, we will be generalizing the Sieve of Eratosthenes. We will be introducing a good amount of notation. For clarity, we'll write down what these correspond to in the Original Sieve of Eratosthenes in {\blue{blue}}.

  Suppose we have a certain subset $A$ of the natural numbers (\blue{ $A = [1, \ldots, x]$ in $\pi(x,z)$}), and $P$ is a set of primes that we're going to use in our sieve (\blue{ $P = P(z)$}). So far, I think these are totally natural thought paths. For each prime in $P$, suppose we have a distinguished set of residue classes. When I mean distinguished set, I just mean that there are some residue classes that we want to sieve out by (\blue{ We cared about getting rid of those numbers that were divisible by some $p \leq z$, so for each prime the distinguished residue class in the standard Sieve of Eratosthenes was $ 0 \bmod p$}). This may be a bit confusing now, but we'll do a nontrivial example in a moment that should make this more clear.

  For each prime, let $\omega(p)$ denote the number of distinguished residue classes for that prime (\blue{ $\omega(p) = 1$ for all $p$ in Eratosthenes}), and let $A_p$ denote those elements in $A$ that are in any of the distinguished residue classes for the prime $p$, or rather $A_p = \left\{ a \in A : a \in \text{ distinguished residue class } \bmod p \right\}$ (\blue{ $A_p$ is the set of numbers up to $x$ divisible by $p$, so $|A_p| = \lfloor x/p \rfloor$}).

  Let's take a moment here to examine something. In the Sieve of Eratosthenes, $|A_p| = \left\lfloor \dfrac{x}{p} \right \rfloor = \left\lfloor \dfrac{x}{p} \right \rfloor \omega(p)$. But the floor function is not a nice function - it's not multiplicative and sometimes has erratic behavior. On the other hand, $\left\lfloor \dfrac{x}{p} \right \rfloor \approx \left( \dfrac{x}{p} \right)$, and in fact isn't really more than $1$ off. This is much better behaved than the floor function. All together, $|A_d| = x \left( \dfrac{\omega(p)}{p} \right) + (\text{small error})$.

  In our generalization, we want $|A_p|$ to be roughly equal to a multiplicative function times $\dfrac{\omega(p)}{p}$. This is an assumption under this sieve. Assuming this is true, write $|A_p| = X(x) \dfrac{\omega(p)}{p} + (\text{small error})$ for some function multiplicative $X(x)$. (\blue{ $X(x) = x$}). For $d$ square-free, let $A_d := \bigcap_{p \mid d} a_p$ and $\omega(d) = \prod_{p \mid d} \omega(p)$ (\blue{ So $A_d$ denotes integers divisible by every prime dividing $d$, useful in the inclusion/exclusion argument, and $\omega(d) = 1$}).

  Finally, call $S(A, P) = |A \setminus \bigcup_{p \in P}A_p|$ the number of elements in $A$ that are not in any of the $A_p$, or equivalently not in any distinguished residue class for any prime (\blue{ $S(A, P(z)) = \pi(x, z)$}). Then we are ultimately counting, or rather estimating, $S(A, P)$ - and it feels very ``sieve-like'' in that we have our elements $A$ and we are taking out elements $A_p$ based on their multiplicative properties. 

  \begin{theorem}[Generalized Sieve of Eratosthenes]
    With $A, P, \omega(p), \omega(d), A_p, A_d,$ and $S(A, P)$ as above, and assuming that $|A_d| = x \left( \dfrac{\omega(p)}{p} \right) + (\text{small error})$, we have
    \begin{equation}
      S(A, P) = X(x)\prod_{p \in P}\left(1 - \frac{\omega(p)}{p}\right) + O(\text{error}),
    \end{equation}
    where the error term is beyond the scope of this talk and paper - any analytic number theory book including sieve theory will mention it. One reference would be Iwaniec-Kowalski \cite{iwaniec2004analytic}.
  \end{theorem}

  \begin{proof}[Sketch of Proof]
    Although it might not feel like it, this result is conceptually no different than the earlier Sieve of Eratosthenes we discussed. The key central bit is understanding that since everything is multiplicative, we can still use inclusion/exclusion. Then if we slightly abuse notation and let $P$ denote the product of the primes we're interested in, we have
    \[
      S(A,P) = \sum_{d \mid P}\mu(d)|A_d| = \sum_{d \mid P} \mu(d) \left( X(x)\frac{\omega(d)}{d} + (\text{error}) \right).
    \]
    Ignoring the error terms (which sieve theorists would say is the most important term to pay attention to), we see that 
    \[
      \sum_{d \mid P} \mu(d)X(x)\frac{\omega(p)}{p} = X(x) \sum_{d \mid P} \mu(d) \frac{\omega(d)}{d}.
    \]
    A basic fact from multiplicative number theory is that if $f$ is a multiplicative function, then $\sum_{n \geq 1} f(n) = \prod_p (f(1) + f(p) + f(p^2) + f(p^3) + \ldots)$, which is really just the fact that integers factor uniquely in disguise. As $\mu(d) \dfrac{\omega(d)}{d}$ is multiplicative, we expect the same here (roughly). Then since $\mu(p) = -1$, we get the minus sign in the final answer, and as $\mu(p^2) = 0$ (and all higher powers), we get only $\left( 1 - \dfrac{\omega(p)}{p} \right)$ per prime in the final answer.
    The interested follower should try to do a more careful analysis, actually paying attention to the error terms along the way, and compare with a reference such as Iwaniec-Kowalski \cite{iwaniec2004analytic}. It's much easier to assume the error terms have a multiplicative bound.
  \end{proof}

  Let's immediately hop into an example application: estimating the number of twin primes.

  \begin{theorem}
    The number of primes $p$ such that $p+2$ is also prime is $O\left( \dfrac{x (\log\log x)^2}{(\log x)^2} \right)$.
  \end{theorem}

  \begin{proof}
    Let $A = [1, \ldots, x]$, and let $P = P(z)$, where I'll specify $z$ later, but not including $2$. So $P = \displaystyle \prod_{2 < p \leq z} p$. For each $p$, distinguish
    \[
      \begin{cases} 0 \bmod p \\
        -2 \bmod p
      \end{cases},
    \]
    so that $\omega(p) = 2$ and $\omega(d) = 2^{\nu(d)}$ (where $\nu(d)$ is the number of prime divisors of $d$). (We exclude the prime $2$ so that it makes easy sense to talk about $-2 \bmod p$). Then $A_p$ includes numbers divisible by primes or 2 less than multiples of primes. A big insight here is that if $p, p+2$ is a pair of twin primes greater than $z$, then since neither $p$ nor $p+2$ is a multiple of a smaller prime, and since $p$ is not 2 less than a multiple of a smaller prime (as that would mean that $p+2$ is a multiple of a smaller prime), $p$ will be counted in $S(A, P)$. Now, $p+2$ might not, and many other things might get in that shouldn't. But since every lower twin prime greater than $z$ is in $S(A, P)$, we have that the number of twin primes in $[z, x]$ is bounded above by $2S(A, P)$. So this sum counts what we want.

    Then $|A_p| = x \dfrac{\omega(p)}{p} + (\text{error})$, as roughly $2$ of every $p$ numbers will be included. Then $|A_d| = x \dfrac{2^{\nu(d)}}{d} + (\text{error})$. 

    Altogether, this means that 
    \[
      S(A, P) = x \prod_{p \mid P(z), p \neq 2}\left( 1 - \frac{2}{p} \right)  + O(\text{error}).
    \]
    As $(1 - \frac{1}{p})^2 = (1 - \frac{2}{p} + \frac{1}{p^2}) > (1 - \frac{2}{p})$, we get that
    \[
      S(A, P) < x \prod_{p \mid P(z), p \neq 2} \left( 1 - \frac{1}{p} \right)^2 + O(\text{error}).
    \]
    It just happens Mertens analyzed the partial product $\prod_{p < z} \left( 1 - \frac{1}{p} \right)^2$ and found that
    \[
      \prod_{p < z} \left( 1 - \frac{1}{p} \right)^2 \approx \left( \frac{e^{-\gamma}}{\log z} \right)^2,
    \]
    where $\gamma$ is the Euler-Mascheroni constant. Using this and a precise knowledge of how the error terms behave in the Sieve of Eratosthenes would lead one to choose $z$ so that $\log z = \dfrac{\log x}{5 \log \log x}$ to minimize the bound while maintaining a sufficiently small error to be meaningful. Notice that this means that $z < x^{1/5}$, and so this shows that the number of twin primes in the range $[x^{1/5}, x]$ is bounded by $S(A, P)$, and
    \[
      S(A, P) = O \left( \frac{x (\log \log x)^2}{(\log x)^2} \right).
    \]

    Even though we omitted the twin primes in $[1, x^{1/5}]$, this doesn't affect the asymptotic. If we assumed that every number in $[1, x^{1/5}]$ was a twin prime, there would be $O \left( \frac{x (\log \log x)^2}{(\log x)^2} \right) + x^{1/5} = O \left( \frac{x (\log \log x)^2}{(\log x)^2} \right)$. And thus if $\pi_2(x)$ represents the number of twin primes up to $x$, then
    \begin{equation}
      \pi_2(x) = O \left( \frac{x (\log \log x)^2}{(\log x)^2} \right).
    \end{equation}

    This concludes the proof.
  \end{proof}

  Notice that this proof didn't really ever rely on the exact residue classes except to identify twin primes. If we were to distinguish $0 \bmod p$ and $-4 \bmod p$ (and omit the prime $3$ as well), this proof would carry through entirely. This gives us the corollary

  \begin{corollary}
    If $\pi_{2n}(x)$ represents the number of pairs of primes $p, p + 2n$ up to $x$, then
    \[
      \pi_{2n}(x) = O \left( \frac{x (\log \log x)^2}{(\log x)^2} \right).
    \]
  \end{corollary}

  We can use this result to prove Brun's Theorem (although in a very different way than Brun proved it himself).

  \begin{theorem}
    \[
      \sum_{p \text{ twin}} \frac{1}{p} < \infty.
    \]
  \end{theorem}

  \begin{proof}
    If I were to summarize this proof with a single phrase, it would be ``partial summation.''

    For those unfamiliar with partial summation (sometimes also called summation by parts), it is integration by parts with Riemann-Stieltjes integrals as opposed to normal Riemann integrals; or alternatively it is integration by parts with more general measures than the typical Lebesgue/Euclidean measure.

    For two relatively well-behaved functions $f, g$, we can think of $\displaystyle \int_a^b f \mathrm d g$ as the limit of the sums $\sum f(x_i) (g(x_{i + 1}) - g(x_i))$, which gives a sort of weight to $f$ based on how quickly $g$ is changing. If $g$ is changing rapidly, then those values of $f$ matter a lot. If $g$ is constant, then the integral is $0$. Notice that when $g(x) \equiv x$, this is exactly a Riemann integral.

    It turns out that integration still works perfectly well with all the normal bells and whistles. Most importantly to us, we can still use integration by parts. With this in mind, we have
    \begin{align*}
      \sum_{p \text{ twin}} \frac{1}{p} &= \sum_{n \geq 1} \frac{1}{n} (\pi_2(n+1) - \pi_2(n)) = \int_1^\infty \frac{1}{\lfloor x \rfloor} \mathrm{d}(\pi_2(x)) = \\
      &= \left[ \frac{\pi_2(x)}{x} \right]_1^\infty - \int_1^\infty \pi_2(x) \mathrm{d}\left( \frac{1}{\lfloor x \rfloor} \right),
    \end{align*}
    and the first term is $0$ because on the one hand, $\pi_2(1) = 0$, and on the other hand we know that $\dfrac{\pi_2(x)}{x} \to 0$ from the asymptotic we proved earlier. We were able to turn from sum to integral's measure changes at discrete steps.

    So we have
    \[
      -\int_1^\infty \pi_2(x) \mathrm{d}\left( \frac{1}{\lfloor x \rfloor} \right) = -\sum_{n \geq 1} \pi_2(n) \left( \frac{1}{n+1} - \frac{1}{n} \right) = \sum_{n \geq 1} \frac{\pi_2}{n^2},
    \]
    where we changed from integral to sum by writing down how the integral changes at those discrete points where it takes value.

    Notice that $\dfrac{\pi_2(n)}{n^2}$ is a positive, decreasing function with limit $0$. Thus we can apply the first-year calculus integral test of convergence (two things for my math 100 students!) to see that this sum converges if and only if the integral $\displaystyle \int_1^\infty \frac{\pi_2(t)}{t^2} \mathrm{d}t$ converges. But
    \[ 
      \int_1^\infty \frac{\pi_2(t)}{t^2} \mathrm{d}t \ll \int_1^\infty \frac{(\log \log t)^2}{t (\log t)^2} \mathrm{d}t \ll \int_1^\infty \frac{\mathrm{d}t}{t (\log t)^{1.5}}  < \infty,
    \]
    by standard $u$-substitution. Thus $\displaystyle \sum_{p \text{ twin}} \frac{1}{p} < \infty$.
  \end{proof}

  Brun successfully proved this theorem, and since the sum is finite, we call the value of the sum Brun's constant. The convergence of this sum is \emph{extremely} slow, so estimates of it are relatively poor and conjectural. But we know that it's approximately $1.9$ - far less then the infinity that some might have hoped for.

  \begin{quote}
    There is an interesting and amusing story unifying some of the constant's we have seen here today. In 2011, Google was bidding on the acquisition of a large set of patents from Nortel. Google's first real bid was \$1,902,160,540, the first 10 digits of Brun's constant. When outbid, Google upped their bid to \$2,614,972,128, the first 10 digits of Merten's second constant (earlier we mentioned Merten's theorem, which is very closely related). When outbid again, Google upped their bid to \$3,141,592,653, the first 10 digits of $\pi$. It turns out that Google lost the auction (it went for about \$4.5 billion), and they've likely regretted it since.  
  \end{quote}

  So we have shown there are ``few'' twin primes. In fact, we expect that $\pi_2(x) \approx \dfrac{cx}{\log^2(x)}$ for a particular constant $c$, and numerical evidence supports this guess. This means that what we have from the Sieve of Eratosthenes is a gross overestimate - we really let too much in. But getting better estimates is hard. One possible strategy would be to sieve through more of the primes up to $x$ - since we only go up to $x^{1/5}$ or so, many non-primes get though the sieve. But this ruins our error estimates and actually worsens our bound. What we really would like is a lower bound instead of an upper bound - but this is largely beyond the Sieve of Eratosthenes.

  Brun actually developed his own sieve (now called Brun's Sieve) to approach the twin primes problem, and his sieve \emph{can} give lower bounds. You might wonder why Brun developed his own sieve rather than using the much older and established Sieve of Eratosthenes. The answer lies hidden in the error analysis that we've omitted from this discussion. Some of the most technical parts rely on results that are younger than Brun. In fact, almost no one touched sieves between Eratosthenes and Brun - there was no interest. But Brun managed to breathe life into the field by giving it new ideas and promise.

  With Brun's Sieve, one \emph{can} prove that there are infinitely many pairs $p, p+2$ where $p$ is prime and $p+2 = P_{20}$, where I use $P_n$ to mean that it's a number with at most $20$ factors (with multiplicity). This is a supposed to be a relatively simple exercise (though I haven't done it, so I can't really say) - it is in Iwaniec-Kowalski. With some improvements, one can prove Chen's Theorem: there are infinitely many pairs $p, p+2$ where $p$ is prime and $p+2 = P_2$ an almost-prime. For a long time, this was the closest anyone got to the twin primes conjecture. (Interestingly, with slight modification, it was also the closest anyone got to the Goldbach conjecture, which are very similar through the lens of Sieve Theory). However it seems unlikely that Chen's results can be improved without significant modification - twin primes are beyond Brun's Sieve.

  The idea behind Brun's Sieve is in essence the same as the Sieve of Eratosthenes, except that one splits apart the sums into positive and negative parts and uses a couple more nice multiplicative function tricks. These give more places for optimization, which is a big part of sieve theory: optimize optimize optimize.

  \section{Modern Results}

  After Brun re\"{e}stablished interest in sieve theory, different and more powerful sieves emerged. Although there are many (google can tell you as much), the one we're going to talk about is Selberg's Sieve. Here's an idea of how we might try to apply Selberg's Sieve to twin primes (and in the process, we'll get an idea of what Selberg's Sieve is). Let
  \[
    \Theta(n) = \begin{cases} \log n & n \text{ prime} \\
      0 & \text{else} \end{cases},
  \]
  and consider the pair of functions
  \[ 
    S_1(x) = \sum_{x < n < 2x} f(n), \qquad S_2(x) = \sum_{x < n < 2x} (\Theta(n+2) - \Theta(n))f(n),
  \]
  for some to-be-chosen function $f(n) \geq 0$. Notice that we are summing across $n$ between some $x$ and $2x$. If both $n$ and $n+2$ are prime, then $\Theta(n) + \Theta(n+2) \approx \log x + \log x = \log x^2$, and in particular 
  \[
    (\Theta(n+2) - \Theta(n))f(n) - \log(3x)f(n) > 0.
  \]
  If we don't have that both $n$ and $n+2$ are prime, then $\Theta(n) + \Theta(n+2) \approx \log x$ (or $= 0$), and 
  \[
    (\Theta(n+2) - \Theta(n))f(n) - \log(3x)f(n) < 0.
  \]

  So we might try to find a function $f(n)$ such that $S_2(x) - \log 3x S_1(x) > 0$ for sufficiently large $x$ or at least for infinitely many choices of $x$. If this is the case, then by the cases mentioned above there must be a pair of twin primes in $[x, 2x]$. In a sense, we have thrown in an additional function $f$ for greater control and optimization - although the path to choosing such $f$ is not at all clear (that's probably what makes it exciting).

  It shouldn't come as a great surprise that we haven't found such an $f$. So we try to get a weaker result. Let $\mathcal{H} = \left\{ h_1, h_2, \ldots, h_k \right\}$ be a set of numbers with a property called ``admissibility'' (i.e. $\mathcal{H}$ is admissible). What this means is that nothing trivial is preventing the numbers $p, p+h_1, p+h_2, p+h_3, \ldots, p + h_k$ from all being prime numbers infinitely often. For example, we will never find a pair of primes $p, p+7$, since they are an odd number apart and the only even prime is $2$. Slightly more meaningfully, we won't find a trio of primes $p, p+2, p+4$ above $3, 5, 7$ because all residue classes $\bmod 3$ are represented, so one of the three numbers will always be divisible by $3$.

  If $\mathcal{H}$ fails one or more of these residue tests, we call $\mathcal{H}$ inadmissible. Goldston (and later Pintz, Yildirim, Zhang, Maynard, and polymath8) considered the pair of functions
  \[
    S_1(x) = \sum_{x < n < 2x} f(n), \qquad S_2(x) = \sum_{x < n < 2x} \left( \sum_{h \in \mathcal{H}}\Theta(n + h) \right) f(n),
  \]
  for $\mathcal{H}$ some admissible set. Then if one could find an $f$ where $S_2 - \log 3x S_1 > 0$ infinitely often, we would have infinitely many prime pairs of the form $p, p+h$ for some $h$ in $\mathcal{H}$. So one might try different functions $f$, bound $S_2$ from below and $S_1$ from above, and see what happens. This more or less happened, even, but saying it now would omit an important part of the story.

  Goldston, Pintz, and Yildirim used this style of sieve, and many very subtle arguments, to prove some results towards the twin prime conjecture (\cite{goldston2006primes}, \cite{goldston2007primes}, \cite{goldston2006small}). Their sieve is often called the GPY sieve after their names. Something of note is that they sieved all the way up to $x$ (not up to $x^{1/5}$ like we did), which means they have really tight control and understanding of error terms. They showed that a relevant piece of information is the ``level of distribution of the primes'', which is a measure of how much the distribution of the primes matches some naive estimates. In particular, if $\forall \epsilon$ there is an $\epsilon' > 0$ such that
  \[
    \sum_{q < x^{\vartheta - \epsilon}} \max_a \left| \pi(x; q,a) - \frac{\text{li}(x)}{\varphi(q)}\right| \ll x^{1 - \epsilon'},
  \]
  where $\text{li}(x) = \int_1^x \frac{\log x}{x} \mathrm{d}x$ is the logarithmic integral, which happens to give an asypmtotically close estimate to the number of primes up to $x$, $\varphi(q)$ is the number of integers less than $q$ that are relatively prime to $q$, and $\pi(x; q,a)$ is the number of primes less than $x$ that are in the arithmetic progression $(a, a + q, a + 2q, \ldots)$, then we say the ``level of distribution of the primes'' is at least $\vartheta$. Dirichlet proved that $\pi(x;q,a) \approx \dfrac{\text{li}(x)}{\varphi(q)}$, so in a sense this is bounding how far away the primes are from our estimate.

  In 1965, Bombieri and Vinogradov showed that $\vartheta \geq \frac{1}{2}$. In 2005, Goldston, Pintz, and Yildirim showed that if $\vartheta > \frac{1}{2}$, then $\liminf(p_n - p_n) < \infty$, where $p_n$ is the $n$th prime. In other words, if $\vartheta > \frac{1}{2}$, then there are infinitely many primes pairs $p, p+h$ for some finite and fixed $h$. They also showed that unconditionally, we're not far off:
  \[
    \liminf_{n \to \infty} \frac{p_{n+1} - p_n}{\log n} \to 0.
  \]
  And this was the direction of progress. Unfortunately, no one knew how to prove anything stronger about $\vartheta$. Some could prove other things - if $\vartheta > 0.971$, then there are infinitely many prime pairs $p, p+h$ for some finite and fixed $h \leq 16$, for example. But there was no other serious avenue of progress.

  Then came Yitang Zhang \cite{zhangbounded}, using a very similar sieve to Goldston's Selberg-style sieve. One way of stating what Zhang did is that he sieved over fewer integers (not all the way up to $x$), but managed to prove new and improved bounds on particularly nasty yet ubiquitious sums called Kloosterman sums, that ultimately allowed him to prevail. A slightly different (and slightly loose) view is that he sieved less efficiently to give him more flexibility, and he controlled the error better than anyone had before.

  Ultimately Zhang proved that $\liminf (p_{n+1} - p_n) \leq 7 \cdot 10^7$ in 2013. His paper is remarkably clear, modular, and easy to read (for a mathematician, that is). He also very openly stated that he did not optimize his result.

  So when Terry Tao started polymath8 \cite{polymath8} with the goal of optimizing Zhang's work, there was rapid progress. They quickly brought $7\cdot10^7$ down to $5414$, and it has continued to decrease since.

  ADDENDUM: I gave this talk before Maynard \cite{maynard2013small} announced his result and before polymath8b started. But it should be mentioned that Maynard (independently) proved this and stronger results, not only about prime pairs, but about arbitrarily large sets of primes occurring infinitely often. Terry Tao's polymath8b started to improve and optimize these new results, and is in the process of doing that \emph{right now}.

  \section{Concluding Remarks}

  I hope this was an enjoyable presentation. There is a list of references at the end containing places for additional information. Many more references can be found at the main site for they polymath projects, which I always endorse and encourage (and sometimes participate in). This note can be found online at \href{http://davidlowryduda.com}{davidlowryduda.com} and at the \href{http://arxiv.org}{arxiv}. 

  \bibliography{A_Friendly_Intro_bibfile}
  \bibliographystyle{plain}

  \end{document}